\documentclass[11pt]{article}
\usepackage{amssymb}
\usepackage{epsfig}
\usepackage{epstopdf}
\usepackage[usenames]{color}
   
   \csname @addtoreset\endcsname{equation}{section}

\usepackage{amssymb}

    \csname @addtoreset\endcsname{figure}{section}
    \csname @addtoreset\endcsname{table}{section}

\newtheorem{theorem}{Theorem}[section]

\newtheorem{proposition}{Proposition}[section]
\newtheorem{corollary}{Corollary}[section]

\newtheorem{remark}{Remark}[section]

\def\d{\delta}
\def\O{\Omega}

\def\Y{{\cal Y}}
\def\F{{\cal F}}
\def\w{\widehat}

\def\Var{{\rm Var\,}}

\def\R{{\bf R}}
\def\E{{\bf E}}
\def\P{{\bf P}}

\def\h{h}

\def\s{\delta}
\def\g{\gamma}

\def\ww{\widetilde}

\def\t{\theta}
\def\oo{\bar}
\def\s{\sigma}

\def\h{h}

\newcommand{\be}{\begin{equation}}
\newcommand{\ee}{\end{equation}}
\newcommand{\bd}{\begin{displaymath}}
\newcommand{\ed}{\end{displaymath}}
\newcommand{\ba}{\begin{array}{ll}}
\newcommand{\ea}{\end{array}}
\newcommand{\baa}{\begin{eqnarray}}
\newcommand{\eaa}{\end{eqnarray}}
\newcommand{\baaa}{\begin{eqnarray*}}
\newcommand{\eaaa}{\end{eqnarray*}}
\font\sm=cmr10

\def\oo{\bar}

\def\mm0{m_{\scriptscriptstyle 0}}
\def\m1{m_{\scriptscriptstyle 1}}

\def\EE{{\mathbb{E}}}

\def\ee{e}
\date{Submitted: June 19, 2015, Revised: April 11, 2017} 
 \title{A pathwise inference method for the parameters of diffusion terms}

\author{
Nikolai Dokuchaev\\ \ {\sm Department of Mathematics \& Statistics,
Curtin University,}\\
{\sm  GPO Box U1987, Perth, 6845 Western Australia}}
\begin{document}
\maketitle
\begin{abstract}   We consider inference
of the parameters of the diffusion term
for  Cox-Ingersoll-Ross  and similar processes with a power type dependence
of the diffusion coefficient from the underlying process.  We suggest some original  pathwise estimates
 for this coefficient and for the power index based on an analysis of  an auxiliary continuous time complex valued process
  generated by the underlying real valued process.
 These estimates  do not rely on the distribution of the underlying process and on a particular choice of the drift.
 Some numerical experiments are used to illustrate the feasibility of the suggested method.
\\
{\bf Key words}: diffusions processes, power index, Cox-Ingersoll-Ross  process, CKLS model, pathwise estimation.
%
\\
{\bf Mathematical Subject Classification (2010):} 91G70       
65C30, 
		65C50, 
		65C60 
\end{abstract}
\def\MM{{\scriptscriptstyle M}}
\section{Introduction}
In this paper, we consider inference  of the diffusion term
for  Cox-Ingersoll-Ross  and similar processes with a power type dependence
of the diffusion coefficient from the underlying process.  These processes are important for applications;
in particular, they are used for interest rate models and for volatility models in finance; see, e.g.,
Heston (1993), Gibbons and Ramaswamy (1993), Lewis (2000), Zhou  (2001),  Carr and  Sun (2007), Andersen and Lund (1997), Gourieroux and Monfort (2013),  Fergusson and Platen (2015), Hin and Dokuchaev (2016),
 and the bibliography therein.  Estimation of the parameters for these models was widely studied; see e.g. Gibbons and Ramaswamy (1993), Andersen and Lund (1997), Kessler  (1997), S{\o}rensen (2000), Zhou  (2001), Fan {\em et al} (2003),  Ait-Sahalia (1996),  De Rossi (2010), Gourieroux and Monfort (2013).

We readdress the problem of inference  for these processes. We suggests  a new method that allows to obtain pathwise estimates of the diffusion coefficient and the power index represented as explicit functions defined on
 an auxiliary continuous time complex valued process generated by the underlying real valued process.  An attractive feature
 of the method is that it not require  to estimate neither
parameters of the drift nor the distributions of the underlying process.  In particular, one does not need to know
the shape of the likelihood function. In addition, our method   allows to consider models with a large number parameters for the drift; therefore, it allows to cover cases where the Maximum Likelihood method is not feasible due to high dimension. This is especially beneficial  for financial application where the trend for the prices is usually  difficult to estimate since it is  overshadowed by  a relatively large volatility. Since the drift is excluded from the analysis,
our method does not lead to an estimation of the drift.  However,
this  could be
  a useful supplement to the existing more comprehensive methods such as described  in
  Gibbons and Ramaswamy (1993), Andersen and Lund (1997), S{\o}rensen (2000), Zhou  (2001), Fan {\em et al} (2003),  Ait-Sahalia (1996),  De Rossi (2010),
Gourieroux and Monfort (2013), and Kessler  (1997). These works used estimation of the parameters for  drift term; on the other  hand, the method discussed in the present paper allows to bypass this task.

Feasibility an robustness of the suggested method is illustrated with some numerical experiments.
\section{The model}
Let  $\t\in\R$ and $T\in (\t,+\infty)$  be given.
We are also given a standard  complete probability space $(\O,\F,\P)$ and a
right-continuous filtration $\{\F_t\}_{t\in[\t,T]}$ of complete $\s$-algebras of
events. In addition, we are given  an one-dimensional Wiener
process $w(t)|_{t\in[\t,T]}$, that is a Wiener
process such that $w(\t)=0$ adapted to  $\{\F_t\}$ and such that  $\F_t$ is independent from $w(s)-w(q)$ if $t\ge s>q\ge \t$.

Consider a  continuous time one-dimensional  random
process $y(t)|_{t\ge \t}$
such that $y(\t)>0$ and \baa  dy(t)=f(y(\cdot),t)dt+\s(t) y(t)^{\g}
d{w}(t),\qquad t\in(\t,T).\label{CT} \eaa
 Here $\g\in [0,1]$,  $\s(t)$  is a bounded  $\F_t$-adapted process,
 $f(s,t):C([\t,T])\times[\t,T]\to\R$
is a measurable function such that $f(s,t)$ is $\F_t$-adapted for any $s\in C([\t,T])$, and
that $f(s_1,t)=f(s_2,t)$ if $(s_1-s_2)|_{[\t,t]}\equiv 0$. In addition, we assume that, for any $\d>0$,
$|f(s_1,t)-f(s_2,t)|\le c_1\|s_1-s_2\|_{C([0,T])}$ and $|f(s,t)|\le c_2(\|s\|_{C([\t,T])}+1)$  for some constants $c_k=c_k(\d)>0$ a.s. (almost surely) a for all $t\in[\t,T]$, $s_1,s_2\in C([\t,T])$, such that $\inf_{t\in[\t,T],k=1,2}s_k(t)>\d$.

Under these assumptions, there exists a Markov time
$\tau$ with respect to $\{\F_t\}$ with values in $(\t,T]$ such that there exists an unique almost surely continuous  solution $y(s)|_{s\in[\t,\tau]}$  such that $\inf_{t\in[\t,\tau]}y(s)>0$.

\subsection*{Examples of applications in financial modelling}
The assumptions on the process $y$ allows to use it for a variety of financial models.
In particular, the assumption on the drift coefficient  $f$ allows to consider a path depending evolution such as described  by equations with delay; see some examples in  Stoica (2005) and  Luong and Dokuchaev (2016).\index{LD,S}

The assumptions on the diffusion coefficient allow to cover many important  financial models. In particular,
 the so-called Cox-Ingersoll-Ross process is used for the interest rate models
and the volatility of stock prices  (Heston (1993)). This equation has the form
\baa dy(t)=a[b-y(t)]dt+\s y(t)^{1/2}
d{w}(t),\qquad t>0,\label{CIR} \eaa  where
 $a>0$, $b>0$, and $\s>0$ are some constants.

A more general model introduced in Chan et al. (1992)
\baa dy(t)=a[b-y(t)]dt+\s y(t)^{\g}
d{w}(t),\qquad t>0,\label{CKLS} \eaa
is called a Chan-Karolyi-Longstaff-Sanders (CKLS) model  in the econometric literature; see e.g. Iacus (2008).
This equation (\ref{CT})  with $\g=2/3$ is also used for volatility modelling;
see, e.g., Carr and  Sun (2007) and  Lewis (2000).

\index{\baaa  dy(t)=\kappa(t)[\mu(t)-y(t)]dt+\s(t)y(t)^{1/2}
d{w}(t)],\quad t>\t.\label{CIR} \eaaa  Here
 $\mu>0$, $\kappa>0$, $\s>0$ are some constants.
It follows from the properties of the solution of (\ref{CT}), as reported in cite{Feller:1951}, that
the following holds:
(1) if $\kappa \equiv 0$ or $\mu\equiv 0$, $y(t)$ reaches zero almost surely and the point zero is absorbing,
(2) if $2\kappa \mu \ge \sigma^2$, $y(t)$ is a transient process that stays positive and never reaches zero,
and
(3) if $0 < 2\kappa \mu < \sigma^2$, $y(t)$ is instantaneously reflective at point zero.}

\section{The main result}\label{SecMain}
Up to the end of the paper, we assume  conditions for  $f$ and $\s$ formulated above holds.
We assume below that $\tau$ is a Markov time with respect  $\{\F_t\}$  such that $\tau\in (\t,T]$ a.s. and
      that $\inf_{s\in[\t,\tau]} y(s)>0$ a.s.. In particular, one can select $\tau= T\land\inf\{s>\t:\  y(s)\le M\}$ for any given
      $M\in (0,y(\t))$.

Our main tool for estimation of the pair $(\s,\g)$ will be provided by the following theorem.
\begin{theorem}\label{ThC} For any $h\in [0,1]$,
 \baa \int^{\tau}_{\t} y(s)^{2(\g-h)} \s(s)^2 ds = 2\log |Y_h(\tau)|\quad\hbox{a.s.}, \label{hY}
\eaa
where $Y_h(s)$ is a complex-valued process defined for $s\in[\t,\tau]$   such that   \baa &&dY_h(s)=iY_h(s)\frac{dy(s)}{y(s)^\h},\quad
s\in(\t,\tau),\nonumber\\&& Y_h(\t)=1.\label{Y}\eaa
\end{theorem}
In (\ref{Y}), $i=\sqrt{-1}$ is the imaginary unit.

\begin{corollary}\label{corr2} (a)  We have that
\baa
\int^{\tau}_{\t} \s(s)^2ds = 2\log |Y_\g(\tau)|\quad\hbox{a.s.}.
\label{sY}
\eaa
 (b) If $\s(t)=\s$ is constant, then, for any $h\in [0,1]$,
 \baa \s^2= 2\left(\int^{\tau}_{\t} y(s)^{2(\g-h)}ds \right)^{-1}  \log |Y_h(\tau)|\quad\hbox{a.s.}. \label{sX}
\eaa
\end{corollary}
\section{Applications of Theorem \ref{ThC} to estimation of $(\g,\s)$}\label{SecAppl}
Up to the end of this paper, we assume that $\s(t)\equiv \s$ is an unknown positive constant.

We present below estimates of $(\s,\g)$  based on available samples $\{y(t_k)\}$, where $t_k\in[\t,\tau\land T]$, such that $t_{k+1}=t_k+\d$ for
$k=\mm0,\mm0+1,...,m-1$,    $\d =(\tau\land T-\t)/(m- \mm0)$, $t_{\mm0}=\t$ and $t_m=\tau\land T$. In this setting,   $y(t_{k})>0$ for $k= \mm0,...,m$.

For $h\in[0,1]$, let  \baaa
 \eta_{h,k}= \frac{y(t_k)-y(t_{k-1})}{y(t_{k-1})^{h}}.
\label{eta}\eaaa

\subsection*{Estimation of $\s$ under the assumption that $\g$ is known}
Let us first suggest an  estimate for $\s$ under the assumption that $\g$ is known.

\begin{corollary}\label{corr3}  For any $h\in[0,1]$, the value $\s$ can be estimated as $\ww\s_{\g,h}$, where  \baa
\ww\s_{\g,h}^2=\left(\d\sum_{k= \mm0+1}^m y(t_k)^{2(\g-h)}\right)^{-1}\sum_{k= \mm0+1}^m \log (1+\eta_{h,k}^2).
\label{ourx} \eaa
\end{corollary}
\subsection*{Estimation of $\g$ with excluded $\s$}\label{SecU}
It appears that Theorem \ref{ThC} implies some useful properties
of the process $Y_h(t)$ allowing to estimate  $\g$ in a setting with unknown constant $\s$.
\begin{proposition}\label{lemmag}  For any  $h_1,h_2\in [0,1]$,
\baa \frac{\int^{\tau}_{\t} y(s)^{2(\g-h_1)}ds}{\int^{\tau}_{\t} y(s)^{2(\g-h_2)}ds}=\frac{\log |Y_{h_1}(\tau)|}{\log |Y_{h_2}(\tau)|}\quad\hbox{a.s.}.
\label{ghh}\eaa
\end{proposition}
\par
Since calculation of  $Y_{h_1}(\tau)$ and $Y_{h_2}(\tau)$ does not require to know the values of $f$, $\g$, and $\s$,
property (\ref{ghh})  allows to calculate
 $\g$ as is shown below.
\begin{corollary}\label{corr4}  An estimate $\w\g$ of $\g$ can be found as a solution of the equation
\baa
\frac{\sum_{k= \mm0+1}^m y(t_k)^{2(\g-h_1)}}{\sum_{k= \mm0+1}^m y(t_k)^{2(\g-h_2)}}= \frac{\sum_{k= \mm0+1}^m \log (1+\eta_{h_1,k}^2)}{\sum_{k= \mm0+1}^m \log (1+\eta_{h_2,k}^2)},
\label{gh} \eaa
for any pair of pre-selected $h_1$ and $h_2$.
\end{corollary}
It can be noted that  $\s$ remains unused and excluded from the analysis for
the method described in Proposition  \ref{lemmag} and Corollary \ref{corr4}; respectively, this method
 does not lead to an estimate of $\s$.
\subsection*{Estimation of the pair $(\s,\g)$}\label{Secg}
\begin{proposition}\label{lemmag2}
The process
\baa \frac{1}{t\land \tau-\t}\log |Y_{\g}(t\land\tau)| \label{cY}\eaa
is a.s. constant in $t\in[\t,T]$.
\end{proposition}
\par
Let
\baaa v_{h,k}=\log(1+\eta_{h,j}^2), \quad k=\mm0+1,...,m,
\label{Loov} \eaaa
and
\baaa \oo  v_h=\frac{1}{m-\mm0}\sum_{j= \mm0+1}^m v_{h,k}.
\label{oL} \eaaa
\begin{corollary}\label{corr5} An estimate  of $\g$ can be found as the solution of the optimization problem
 \baa \hbox{Minimize}\quad  \sum_{k=\mm0+1}^m\left(v_{h,k}- \oo v_h\right)^2\quad
 \hbox{over}\quad h\in[0,1].
 \label{L}\eaa
 In this case, $\s$ can be estimated as
\baa
\w\s=\sqrt{\oo v_{\w\g}/\d}.
\label{wvol}
\eaa
where $\w\g$ is the estimate of $\g$ obtained as a solution of (\ref{L}).
\end{corollary}
\begin{remark}\label{corr05} Corollary \ref{corr5} allows the following modefication
 a special case of estimation of $\g$ for the case where $\s$ is a known constant $\s$: an estimate $\w\g$ of $\g$ can be found
as the solution of the optimization problem
 \baaa \hbox{Minimize}\quad  \sum_{k=\mm0}^m\left(v_{h,k}/\d- \s^2\right)^2\quad
 \hbox{over}\quad h\in[0,1].
 \label{L2}\eaaa
\end{remark}

\section{Proofs}
 For $M\in (0,y(\t))$, let $\tau_M=\tau\land\sup\{s\in[\t,T]:\ \inf_{q\in[\t,s]} y(q)\ge M\}$.
 Clearly,
 \baa \tau_M\to \tau\quad\hbox{as}\quad M\to 0\quad\hbox{a.s.}
 \label{lim}\eaa


\par
{\em Proof of Theorem \ref{ThC}.} The proof follows the idea of the proof of
Lemma 3.2 from Dokuchaev (2014), where  less general log-normal underlying processes were considered. Let $\ww a(t)=f(y(t),t)y(t)^{-h}$.
We have, for any  $M\in (0,y(\t))$,
 \baaa  dY_h(t)=iY_h(t)[\ww a(t)dt+y(t)^{\g-h}\s(t)
d{w}(t)], \quad t\in(\t,\tau_M).\label{CT22} \eaaa
By the Ito formula again,
for any $M\in (0,y(\t))$,  \baaa Y_h(\tau_M)&=&Y_h(\t)\exp\left(i\int_{\t}^{\tau_M}\ww a(s)ds-\frac{i^2}{2}\int^{\tau_M}_{\t}y(s)^{2(\g-h)}\s(s)^2ds
+i\int^{\tau_M}_{\t}y(s)^{\g-h}\s(s)dw(s)\right)\nonumber\\&=&\exp\left(i\int_{\t}^{\tau_M}\ww a(s)ds+\frac{1}{2}\int^{\tau_M}_{\t}y(s)^{2(\g-h)}\s(s)^2ds
+i\int^{\tau_M}_{\t}y(s)^{\g-h}\s(s)dw(s)\right)\quad\hbox{a.s.}. \label{solSx} \eaaa Hence
\baaa |Y_h(\tau_M)|=\exp\left(\frac{1}{2}\int^{\tau_M}_{\t}y(s)^{2(\g-h)}\s(s)^2ds
\right) \quad\hbox{a.s.}\label{SYxx} \eaaa and \baaa\int^{\tau_M}_{\t}y(s)^{2(\g-h)}\s(s)^2ds =2\log |Y_h(\tau_M)|\quad\hbox{a.s.}. \label{sX1}\eaaa
Hence (\ref{hY}) follows from (\ref{lim}). $\Box$

{\em Proof of Corollary \ref{corr2}} follows immediately from Theorem \ref{ThC}. $\Box$
\par
{\em Proof of Corollary \ref{corr3}}.  Let  $t_m=\tau_M$, $t_{ \mm0}=\t$,  and let $t_k=t_{\mm0}+(k-\mm0)\d$ if $\mm0\le k\le m$.  Let  $\eta_{h,k}$ be defined by
(\ref{eta}). The Euler-Maruyama  time  discretization   of (\ref{Y}) leads to the stochastic difference equation
 \def\YY{{\mathcal{Y}}}
 \baaa
&&\YY_h(t_k)=\YY_h(t_{k-1})+i\YY_h(t_{k-1})\eta_{h,k},\quad k\ge \mm0+1, \nonumber\\&&
\YY_h(t_{ \mm0})=1. \label{Xd}\eaaa
  (See,  e.g.,  Kloeden and Platen (1992), Ch. 9). This equation can be rewritten as
 \baaa
&&\YY_h(t_k)=\YY_h(t_{k-1})(1+i\eta_{h,k}),\quad k\ge \mm0+1,\nonumber \\&& \YY_h(t_{
\mm0})=1. \label{Xdd}\eaaa Hence \baaa \YY_h(t_m)=\prod_{k= \mm0+1}^m
(1+i\eta_{h,k}). \eaaa Clearly, \baaa |\YY_h(t_m)|=\prod_{k= \mm0+1}^m
|1+i\eta_{h,k}|=\prod_{k= \mm0+1}^m (1+\eta_{h,k}^2)^{1/2}, \eaaa and \baa
\log|\YY_h(t_m)|=\sum_{k= \mm0+1}^m \log
[(1+\eta_{h,k}^2)^{1/2}]=\frac{1}{2}\sum_{k= \mm0+1}^m \log
(1+\eta_{h,k}^2). \label{Y=}\eaa Then (\ref{sY}) leads to
 estimate (\ref{ourx}). $\Box$

 {\em Proof of Proposition \ref{lemmag} and Proposition \ref{lemmag2}} follows immediately from Theorem \ref{ThC} and Proposition \ref{corr2}(b). $\Box$

 {\em Proof of Corollary \ref{corr4}} follows from the natural discretization of integration and (\ref{Y=}). $\Box$

{\em Proof of Corollary \ref{corr5}}.
 It
follows from (\ref{Y=}) that the sequence $\{\log|\Y_{h}(t_k)|\}$ represents the discretization of the continuous time process
$\log|Y_h(t\land\tau)|$ at points $t=t_k$; this process is linear in time for $h=\g$ and
$2\log|Y_\g(t\land\tau)|\equiv (t\land\tau-\t)\s^2$.
Hence \baaa
2\log|Y_\g(t_{k+1})|-2\log|Y_\g(t_k)|=\d\s^2.\eaaa
On the other hand, (\ref{Y=}) implies that   \baaa
2\log|\Y_h(t_{k+1})|-2\log|\Y_h(t_k)|=v_{h,k},\quad h\in[0,1].
\eaaa
This leads to an optimization problem
\baaa   \hbox{Minimize}\quad \sum_{k=\mm0+1}^m (v_{h,k}/\d- c)^2
\quad \hbox{over}\quad h\in[0,1],\quad c>0.
 \label{Lu}\eaaa
 By the properties of quadratic optimization, this problem
 can be replaced by the problem
 \baa  \hbox{Minimize}\quad \sum_{k=\mm0+1}^m (v_{h,k}/\d- \oo v_h)^2
 \quad \hbox{over}\quad h\in(0,1].
 \label{Lp2}\eaa
Then the proof follows. $\Box$

{\em Proof of Remark \ref{corr05}} repeats the previous proof without optimization over $c$.  $\Box$

\section{Numerical experiments}
To illustrate numerical implementation of the algorithms described above,
we applied these algorithms for discretized Monte-Carlo simulations of some generalized version the Cox-Ingersoll-Ross process (\ref{CIR}). We consider a toy example of a process with a large number of parameter. Presumably,  estimation of all these
parameters is not feasible due a high dimension for a method of moments or Maximum Likelihood Method.

We consider a process evolving as the following:
\baa dy(t)=H\left(y(t),y(\max(t-\lambda,0)\right)dt
+\s y(t)^{\g}
d{w}(t),\qquad t>0,\label{CIRn} \eaa
where
\baaa
&&H(x,y)=\sum_{k=1}^N\left[F_k(x)+G_k(y)\right],\\
&&F_k(x)=a_k[b_k-x^{\nu_k+1/2}]+c_k\cos(d_k x+e_k), \quad  G_k(x)=0.1\,\w a_k[\w b_k-x^{\w\nu_k+1/2}],
\eaaa
The parameters $N,a_k,b_k,\nu_k,c_k,d_k,e_k,\w a_k,\w b_k,\w e_k,\nu_k,\w\nu_k,\lambda$  are
randomly selected in each experiment.
In particular, the integers  $N$ are selected randomly at the set $\{1,2,3,4,5\}$ with equal probability.
The delay parameter $\lambda$ has the uniform distribution on the interval $[0,0.2]$.
The parameters  $a_k,b_k,\nu_k,c_k,d_k,e_k,\w a_k,\w b_k,\w e_k,\nu_k,\w\nu_k$  are
uniformly distributed on the interval $[0,1]$.

For the Monte-Carlo simulation, we considered corresponding discrete time process $\{y(t_k)\}$ evolving as \baa y(t_{k+1})=y(t_{k})+H(y(t_k), y(t_{\max(k-\ell,0)})\d +\s
y(t_k)^{\g}\d^{1/2}\,\xi_{k+1}, \quad k=0,...,n,\label{myy}\eaa
with mutually independent random variables
$\xi_k$ from the standard normal distribution $N(0,1)$. Here
 $\d=t_{k+1}-t_k=1/n$
 this corresponds to $[\t,T]=[0,1]$ for
continuous time underlying model.   The delay $\ell$ is the integer part of $\lambda (T-\t)/(n+1)$.

We considered  $n\in\{52,250, 10000,20000\}$. For the financial applications, the choice
of $n=52$ corresponds to weekly sampling; the choice
of $n=250$ corresponds to daily sampling.

In the Monte-Carlo simulation trials, we considered  random $y(t_{\mm0})$ uniformly distributed on $[0.1,10]$ and truncated paths $y(t_k)|_{\mm0\le k\le m}$, with the Markov stopping time
$m=n\land \inf\{k:\ y(t_k)\le 0.001y(t_{\mm0})\}$.  In this case,  $y(t_{k})>0$ for $k= \mm0,...,m$. To exclude the possibility that $y(t_m)\le 0$ (which may happen for our discrete time process since
the values of $\xi_k$ are unbounded), we replace  $y(t_m)$ defined by (\ref{myy}) by  $y(t_m)=y(t_{m-1})>0$ every time when $m<n$ occurs. It can be noted that, for our choice of parameters, the occurrences  of the event $m<n$ were very rare and have not an impact on the statistics.

We used 10,000  Monte-Carlo trials for each trial (i.e. for each entry in each of the tables \ref{tab1}-\ref{tab3} below).
We found that enlarging the sample
does not improve the results. Actually, the experiments with
 5,000 trials or even 1,000 Monte-Carlo trials produced the same results.

The parameters of the errors obtained in these experiments are
quite robust with respect to the change of other parameters as well.

We denote by  $\EE$
the sample means of the  corresponding values over
all Monte-Carlo simulation trials. For the estimates $(\w\s,\w\g)$ of $(\s,\g)$,
we evaluated  the root mean-squared errors  (RMSE) $\sqrt{\EE \left|\w\s-\s\right|^2}$ and
 $\sqrt{\EE \left|\w\g-\g\right|^2}$,
 the mean errors $\EE |\w\s-\s|$ and $\EE \left|\w\g-\g\right|$,
  and the biases $\EE (\w\s-\s)$ and $\EE \left(\w\g-\g\right)$.

  In the experiment described below, we used $\s=0.3$, $\g=1/2$, and $\gamma=0.6$.

\index{\baaa  \d=1/52,\quad  \d=1/250,\quad \d=1/6000.\label{as11}\eaaa}
\subsection*{Estimation of $\s$ using Corollary \ref{corr3}}
The numerical implementation of Corollary \ref{corr3}
requires  to use the value $\g$.   In other words,  one
have to use  certain  hypothesis about
the value of $\g$, for instance, based on estimation of $\g$ that was done separately.
This setting leads to an error caused by miscalculation of $\g$.

To illustrate the dependence of the error for the estimate of $\s$ from the error in the hypothesis on $\g$,  we considered 
estimates for inputs simulated with $\g=1/2$ and with
different $h$.

 Tables \ref{tab1} (a),(b) show the parameters of the errors in the
experiments described
 above  for estimate (\ref{ourx}) with different $h$ and with $\g=1/2$, for $\d=1/52$ and $\d=1/250$ respectively.
 Numerical experiments shows that these estimates are robust with respect to small  errors for $\g$;
 however, the estimation error for  $\s$  caused by misidentification of $\g$ can be significant.

\subsection*{Estimation of $\g$ with unknown $\s$ using (\ref{gh}) and (\ref{L})}
In these experiments, we used  simulated process with $\g=0.6$ and estimates (\ref{gh}) and (\ref{L}).

For solution of  equation  (\ref{gh}) and optimization problem (\ref{L}),
we used simple search  over a finite set $\{h_k\}_{k=1}^N=\{k/N\}_{k=1}^N$. We used  $N=300$ for (\ref{gh}) and $N=30$ for (\ref{L}). Further increasing of  $N$ does not improve the results but slows down calculation.

It appears that estimation of $\g$ is more numerically challenging than estimation of $\s$
using (\ref{ourx}) with known $\g$. In our experiments,
we observed that the dependence of the value of criterion function in (\ref{L})
 depends on $h$ smoothly and the dependence on $h$ for each particular Monte-Carlo trial is represented by an U-shaped smooth convex function. However, the minimum point of this functions is deviating significantly for different Monte-Carlo trials, especially in the case of low-frequency  sampling.  It requires high-frequency sampling to reduce the error $\w\g-\g$.   Table \ref{tab2} shows the parameters of the error  $\w\g-\g$.
We found that these  parameters are
quite robust with respect to the change of other parameters of simulated process.

\subsubsection*{Estimation of $\s$  using (\ref{L})}
The solution of optimization problem (\ref{L}) gives an estimate of $\s$, in addition to an estimate of $\g$, in the setting with unknown $\s$, via
(\ref{wvol}). This gives a method for estimation of $\s$ tat can be an alternative  to estimator (\ref{ourx}). Table  \ref{tab3} shows the  parameters of the error $\w \s-\s$. It appears that the RMSE is
larger than for estimators (\ref{ourx}) applied with a correct  $h=\g$ and has the same order as the RMSE for this estimators applied with $h\neq \g$, i.e., if $\g$ is "miscalculated".

\index{It can be added that the error of the same estimator (\ref{wvol}) for $\s$  is much smaller if this estimator is  applied with fixed $h=\g$.  In some experiments not detailed here \index{(>> [gL]=gammL(10000,1,10000,30,0)WITH KNOWN GAMMA)} we found that RMSE  for estimator (\ref{wvol}) applied with $h=\g$  is close  to the RMSEs for (\ref{ourx})   with  $h=\g$.}

\begin{table}[ht]
\vspace{-0.cm}\begin{center}(a) $\d=1/52$\\\begin{tabular}
 {|c|c|c|c|}\hline
 $\hphantom{xxxx}$   &$\hphantom{\Biggl|}$ $\sqrt{\EE \left|\w\g-\g\right|^2}$ &$\EE \left|\w\g-\g\right|$ &$\EE \left(\w\g-\g\right)$
 \\
    \hline  $\g=0.5$, $h=0.5$  & 0.0312 &  0.0248 & 0.0034 \\
    \hline  $\g=0.4$, $h=0.5$  &   0.0458  &    0.0365 & 0.0281 \\
\hline  $\g=0.6$, $h=0.5$  &   0.0358  & 0.0290  &  -0.0183 \\
\hline  $\g=0.7$, $h=0.5$  &  0.0495 &  0.0413  &  -0.0370\\
\hline
 \end{tabular}\end{center}
 \begin{center}(b) $\d=1/250$\\\begin{tabular}
 {|c|c|c|c|}\hline
 $\hphantom{xxxx}$   &$\hphantom{\Biggl|}\sqrt{\EE \left|\w\g-\g\right|^2}$ &$\EE \left|\w\g-\g\right|$ &$\EE \left(\w\g-\g\right)$
 \\
\hline  $\g=0.5$, $h=0.5$  &0.0136  &  0.0109 & 0.0006  \\
   \hline  $\g=0.4$, $h=0.5$ & 0.0328 & 0.0272 &  0.0259 \\
\hline  $\g=0.6$, $h=0.5$  &0.0269  &   0.0227  & -0.0215
\\
\hline  $\g=0.7$, $h=0.5$  &  0.0468  &  0.0416 &  -0.0414\\
\hline
 \end{tabular}
\end{center}\vspace{0mm}
\vspace{-0mm}\caption{Parameters of the error $\w\s-\s$  for $\w\s$ obtained from estimates  (\ref{ourx})  with $\d=1/52$ and $\d=1/250$. In the first column, $\gamma$ is the "true" power  used for simulation, and $h$
is the parameter  of (\ref{ourx}) used for estimation; mismatching of $\g$ and $h$ leads to a larger bias and a larger estimation error.  }
  \label{tab1}
  \end{table}

    \begin{table}[ht]
\vspace{-0.cm}\begin{center}\begin{tabular}
 {|c|c|c|c|c|c|c|}\hline
 $\hphantom{xxxx}$   &$\hphantom{\w{\Biggl|}}$ $\sqrt{\EE \left|\w\g-\g\right|^2}$& $\EE |\w\g-\g|$ &$\EE \left(\w\g-\g\right)$ &$\sqrt{\EE \left|\w\g-\g\right|^2}$ &$\EE \left|\w\g-\g\right|$ &$\EE (\w\g-\g)$\\
& for (\ref{gh}) & for (\ref{gh})& for (\ref{gh}) & for  (\ref{L}) &  for  (\ref{L}) & for  (\ref{L})\\
   \hline
 $\d=1/250$     &  0.2078 &0.1736 & 0.1078 & 0.2304 &0.1946 & 0.1166 \\
   \hline
  $\d=1/10,000$     &  0.0309 &0.0182 & 0.0039 & 0.0356 &0.0221 & 0.0042 \\
   \hline $\d=1/20,000$     &  0.0222 &0.0109 & 0.0020 & 0.0483 &0.0294 & 0.0004 \\
   \hline

 \end{tabular}
\end{center}\vspace{0mm}
\caption{Parameters of the error $\w\g-\g$  for the solution of (\ref{gh})
and   (\ref{L}) with an unknown $\s$.}
  \label{tab2}
  \end{table}
  \begin{table}[ht]

\vspace{-0.cm}
\begin{center}\begin{tabular}
 {|c|c|c|c|}\hline
 $\hphantom{xxxx}$   &$\hphantom{\Biggl|}\sqrt{\EE \left|\w\s-\s\right|^2}$ &$\EE \left|\w\s-\s\right|$ &$\EE (\w\s-\s)$\\
   \hline
 $\d=1/250$ &  0.0515& 0.0264 &0.0092\\
   \hline
  $\d=1/10,000$ & 0.0063&  0.0038&  0.0001\\  \hline
  $\d=1/20,000$ & 0.0168&  0.0108&  0.00003 \\
    \hline
 \end{tabular}
\end{center}\vspace{0mm}
\caption{Parameters of the error $\w\s-\s$  for $\w\s$ obtained from (\ref{wvol}) and (\ref{L}) with an unknown $\g$.}
  \label{tab3} \end{table}

\newpage
\subsection*{Comparison with the performance of other metods} S{\o}rensen (2000) \index{, p.94 and p. 98,} and Zhou (2001) reported the results of testing of a variety of  estimators  based on the maximum likelihood method or the method of moments for  special cases of (\ref{CT}). These works 
considered  simulated processes with  a preselected structure for the drift term with a low dimension of the vector of parameters. Due to numerical challenges for the methods used, the number of Monte Carlo trials was relatively short in these works 
(100 trials in S{\o}rensen (2000) and 1,000 trials in Zhou (2001)).
 S{\o}rensen (2000) considered model (\ref{CKLS})  with one fixed set of parameters $(a,b)$ for the drift, 
 and Zhou (2001) considered  model  (\ref{CIR}) for a variety of  the  parameters  $(a,b)$ for the drift. 
 S{\o}rensen (2000) considered estimation of  $(\s,\gamma)$ and estimation of the drift parameters,  and Zhou (2001) considered estimation of $\s$ and estimation  of the drift parameters with fixed $\g=1/2$.

 The results for $\s$  in Table 5 from  Zhou (2001) reported for $\d=1/500$ depends significantly on the choice of the 
 the drift parameters $(a,b)$ in (\ref{CIR}) (in our notations). The minimal RMSE  for estimates of $\s$ among all pairs $(a,b)$ is of the same order as the RMSE  reported in our Table \ref{tab1}(a) for $\d=1/250$   for the case of known $h$; for other choices of the drift the RMSE in Table 5 from  Zhou (2001) is much larger.  Remind that RMSE is smaller for smaller $\d$.
 
The RMSE  for $\s$  reported  in Table II.1 from  S{\o}rensen (2000) for $\d=1/500$ (in our notations)   is approximately the same as in Table \ref{tab3}
for $\d=1/250$.   However,   the RMSE for $\g$ with $\d=1/500$ is three times smaller in  Table II.1 from  S{\o}rensen (2000) for some estimators than in Table \ref{tab2} with $\d=1/250$. However, it may happen that
the performance of the estimators  in Table II.1 from  S{\o}rensen (2000) is not robust with respect to different choices  of the drift parameters, similarly to the case presented in  Table 5 from  Zhou (2001) for $\g=1/2$. On the other hand, our method  allowed to include  
a high variety of drift models with almost unlimited dimension, and, as we found in some unreported experiments, the choice of particular drifts 
does not have an impact on the performance of the estimator.  
 
 \section{On the consistency of the method} 
Let us describe briefly the consistency of the method. Clearly, one cannot  the real life data such as market prices are  generated by model (\ref{CT}). Hence we restrict our consideration by the error for simulated data. The equations in the continuous time used for our method  are  exact and hold  almost surely for continuous time underlying processes (\ref{CT}).  Therefore,  the only source of the error  is the time discretization error. This error is inevitable since the the method requires
pathwise evaluation of stochastic integrals. Let us discuss  briefly
 consistency of the method as convergence  of the estimates 
   to the true values  as the sampling frequency is increasing, i.e.  $\d\to 0$. A rigorous  analysis of convergency for the time discretization  requires significant analytical efforts outside of the scope of this paper; see e.g. Kloeden and Platen (1992)   and Jourdain and  Kohatsu-Higa  (2011), where  review of the recent literature can be found. 
   We leave this analysis for the future research and give below a short sketch of two possible approaches. 
   
There are two options. First, one can consider Euler-Maruyama  time  discretization  for the pair $(y,Y_h)$
such as described for the numerical experiments described above.  In this case,  $f$ and the sampling frequency $\d$ have to be such that a satisfactory approximation is achieved.  In particular, by Theorem 9.6.2
from  Kloeden and Platen (1992), p. 324, these conditions are satisfied for CIR models  as well as for the case
where $f(y(\cdot),t)=f(y(t),t)$.  Some  analysis o
 and  conditions for the convergence  in more general cases can be found in Jourdain and  Kohatsu-Higa  (2011).
The numerical experiments described above demonstrate that the required convergence takes place for 
equtations will delay modelled there.

Another option is to consider convergence of  the method for $\d\to0$  given that $\Y(t_k)$ are 
constructed with  the "true" entries $y(t)$.  We presume here that it is possible to produce an arbitrarily close approximation of  a continuous path $y(t)$ via Monte-Carlo simulation with increasing of the simulation frequency. 
Let $t_\d(t)$ be selected as $t_k$ such that $|t-t_k|=\min_p|t-t_p|$ (for certainty, let it me the minimal point if $t$ is in the middle of a sampling interval). 
   Clearly,  $\E \sup_{t\in [0,\tau]}|y(t_\d(t))-y(t)|^2\to 0$  as $\d\to 0$. Hence 
   $\d\sum_{k= \mm0+1}^m y(t_k)^{2(\g-h)}\to \int^{\tau}_{\t} y(s)^{2(\g-h)}ds$  in probability as $\d\to 0$. 
  Further,  can be shown that  $\log|\Y_h(t_\d(t))|\to \log|Y_h(t)|$  in probability as $\d\to 0$.
  This leads to converges of estimates to their true values in probability.  \section{Discussion}
  \begin{enumerate}
\item
The estimates listed  in Section \ref{SecAppl} do not use neither  $f$ nor the probability distribution of the process $\{y(t)\}$. In particular, they are invariant with respect to the choice of an equivalent probability measure.  This is an attractive feature that allows to consider models with a large number of parameters for the drift.
\item It appears that estimation of the power index $\g$ with unknown $\s$ is numerically challenging
 and requires  high-frequency sampling to reduce the error.  Perhaps, this can be improved using other modifications of (\ref{L}) and other estimates for the degree of nonlinearity for the implementation of Lemma \ref{lemmag2}. In particular, the standard criterions of
linearity for the first order regressions could be used, and  $L_2$-type
 criterions could be replaced by $L_p$-type criterions with $p\neq 2$. So far, we were unable to find a way to reduce the error for lower sampling frequency.
We leave it for the future research.
\item
Our approach does not cover the estimation of the drift $f$ which is a more
 challenging problem. However, the estimates for $(\s,\g)$  suggested above can be used to simplify statistical inference  for $f$ by reduction of the dimension of the numerical problems arising in the maximum likelihood method, methods of moments,  or least squares estimators, for $(f,\s,\g)$. This can be illustrated as the following.
\par  Assume that $\g=1/2$ is given and that evolution of  $y(t)$  is described by Cox-Ingersoll-Ross equation (\ref{CIR}) with $\t=0$.
It is known that
\baaa
&&\E y(T)=b(1-e^{-aT})+e^{-aT}y(0),\nonumber\\
&& \Var y(T)=\frac{\s^2}{2a}b(1-e^{-aT})+e^{-aT}\frac{\s^2}{a^2}(1-e^{-aT})y(0).
\label{V}\eaaa
(See, e.g., Gourieroux and Monfort (2013)). This system can be solved with respect to $(a,b)$  given that $\E y(T)$ and $\Var y(T)$ are estimated by their sampling values, and $\s$ is
estimated as suggested  above.
\index{This can be rewritten as
 \baa
&& b(1-e^{-aT})=\E y(T)-e^{-aT}y(0),\label{E}\\
&& \Var y(T)=\frac{\s^2}{2a}\left[\E y(T)-e^{-aT}y(0)\right]+e^{-aT}\frac{\s^2}{a^2}(1-e^{-aT})y(0).
\label{V1}\eaa
}
\item The paper focuses on the case where $\s(t)$ is constant. However,
some results can be extended on the case of time depending and random $\s(t)$. For example, the proofs given  above
 imply  that   $\int^{t_m}_{\t} \s(s)^2ds$ can be estimated as
 \baaa
\ww\s_{\g,\g}^2= \sum_{k= \mm0+1}^m \log (1+\eta_{\g,k}^2).
\label{oury}
\eaaa
 \end{enumerate}
\subsection*{Acknowledgments}  The author gratefully acknowledges support provided
by ARC grant of Australia DP120100928.

\end{document}